\documentclass[a4paper,10pt]{article}
\usepackage{amsthm}
\usepackage{amsmath}
\usepackage[margin=1.2in]{geometry}
\usepackage{amsfonts}
\usepackage{amssymb}
\usepackage{latexsym}
\usepackage{graphicx}

\usepackage{tabularx}
\usepackage{enumitem}
\usepackage{ae}
\usepackage[T1]{fontenc}

\usepackage{geometry}
\newcommand{\Keywords}[1]{\par\noindent{\bfseries Keywords}: #1}
\newcommand{\MSC}[1]{\par\noindent{\bfseries MSC}: #1}

\newcommand{\PaperTitle}[1]{#1}

\newcommand{\JournalName}[1]{\emph{#1}}

\newtheorem{prop}{Proposition}[section]
\newtheorem{lem}[prop]{Lemma}
\newtheorem{cor}[prop]{Corollary}
\newtheorem{thm}[prop]{Theorem}

\newtheorem{pro}[prop]{Problem}
\newtheorem{prob}[prop]{Open Problem}
\theoremstyle{definition}
\newtheorem{rem}[prop]{Remark}
\newtheorem{defin}[prop]{Definition}

\numberwithin{equation}{section}

\date{}
\author{Simone Costa \thanks{DII/DICATAM - Sez. Matematica, Universit\`a degli Studi di Brescia, Via
Branze 38, I-25123 Brescia, Italy, email: simone.costa@unibs.it}
\thanks{This work was partially supported by Italian Ministry of Education under grant PRIN 2015 D72F16000790001}}
\title{A complete solution to the infinite\\ Oberwolfach problem}
\begin{document}
\maketitle
\begin{abstract}
Let $F$ be a $2$-regular graph of order $v$.
The Oberwolfach problem, $OP(F)$, asks for a $2$-factorization of the complete graph on $v$ vertices in which each $2$-factor is isomorphic to $F$.
In this paper, we give a complete solution to the Oberwolfach problem over infinite complete graphs, proving the existence of solutions that are regular under the action of a given involution free group $G$. We will also consider the same problem in the more general contest of graphs $F$ that are spanning subgraphs of an infinite complete graph $\mathbb{K}$ and we provide a solution when $F$ is locally finite. 
Moreover, we characterize the infinite subgraphs $L$ of $F$ such that there exists a solution to $OP(F)$ containing a solution to $OP(L)$.
\end{abstract}
\vskip-0.1cm
\Keywords{ Regular Factorizations, Oberwolfach Problem, Subsystems.}
\MSC{05C70}
\par\vskip.5cm
\section{Introduction}
In this paper we deal with graphs, finite or infinite, which are simple and with no loops. Given a graph $\Lambda$ we denote by $V(\Lambda)$ and $E(\Lambda)$ the set of vertices and the set of edges of $\Lambda$, respectively. As usual we will use the notation $\mathbb{K}_V$ for the complete graph whose vertex set is $V$.

If we define the degree of a vertex $v$ the number of edges $e$ such that $v\in e$, a graph $F$ is $2$-regular if and only if every vertex has degree $2$. 
A spanning subgraph $F$ of a graph $\Lambda$ is a subgraph such that $V(F)=V(\Lambda)$; equivalently a spanning subgraph of $\Lambda$ is a subgraph obtained by edge deletions only. We are now able to provide the definitions of $2$-factor and $2$-factorization:
\begin{defin}
\begin{itemize}
\item A a $2$-factor $F$ of $\Lambda$ is a $2$-regular spanning subgraph of $\Lambda$.
\item A decomposition $\mathcal{C}$ of a graph $\Lambda$ is a partition of the edges of $\Lambda$ into subgraphs. In the case each subgraph is a $2$-factor the decomposition is said to be a $2$-factorization.
\end{itemize}
\end{defin}
Keeping those definitions in mind, we introduce the famous Oberwolfach problem: 
\begin{pro}
Given a $2$-factor $F$ of a complete graph $\mathbb{K}$, the Oberwolfach problem $OP(F)$ asks for a $2$-factorization 
of $\mathbb{K}$ whose $2$-factors are all isomorphic to $F$.
\end{pro}
We remark that a finite $2$-regular graph $F$ is an union of cycles and two finite $2$-regular graphs $F_1$ and $F_2$ are isomorphic if and only if there exists a bijection $\sigma$ between their cycles that preserves their lengths. In this case we may assume $F$ consists of $m_i$ cycles of length $k_i$ for $i=1,\dots,t$ and the $OP(F)$ is denoted by:
$$ OP(k_1^{m_1},\dots, k_t^{m_t}).$$ 

The original formulation, given by Ringel in 1967, asks to arrange a series of meals for an odd number $v$ of people around $t$ tables of sizes $k_1,\dots, k_t$ so that each person sits next to each other exactly once. Despite its simple formulation, a complete solution to the classical $OP$ have not been achieved yet. 

Some important existence results on that problem have been obtained in the so called equipartite case, i.e. the $OP(k^m)$ (Liu and Lick, 2003 \cite{P6}), in the two table case $OP(k_1,k_2)$ (Traetta, 2013 \cite{P7b}) and in the case $|V(F)|$ belongs to an infinite set of prime numbers (Bryant and Scharaschkin, 2009 \cite{P3}) or $|V(F)|$ is finite and big enough (Glock, Joos, Kim, K\"{u}hn, Osthus, 2018 \cite{Asimtp}). It is also known that, up to $4$ exceptions, every instance of the problem has a solution whenever $|V(F)| \leq 60$ (Salassa, Dragotto, Traetta, Buratti and Della Croce, 2019 \cite{SDTBC}; see also Deza, Franek, Hua, Meszka and Rosa, 2010 \cite{P5b}).
A variant of the $OP$ can be considered when $v$ is even: here one can look for $2$-factorizations of the so called cocktail party graph. This case has been solved whenever all cycles have even length (Bryant and Danziger, 2011 \cite{P2}; see also Haggkvist, 1985 \cite{H85}).
In this paper, we consider also the following more general problem:
\begin{pro}
Given a spanning subgraph $F$ of a complete graph $\mathbb{K}$, the Generalized Oberwolfach problem, or briefly Generalized $OP(F)$, asks for a decomposition 
of $\mathbb{K}$ whose members are spanning subgraphs all isomorphic to $F$.
\end{pro}
An interesting related problem is the existence of $2$-factorizations with subsystems. More precisely, let us consider a $2$-regular graph $F$, a solution $\mathfrak{F}$ to the $OP(F)$ over the complete graph $\mathbb{K}_V$, and a $2$-regular subgraph $L$ of $F$. We say that a solution $\mathfrak{L}$ to the $OP(L)$ over $\mathbb{K}_W$ is a subsystem of $\mathfrak{F}$ if $W\subset V$ and 
$\mathfrak{L}$ is the restriction of $\mathfrak{F}$ on $\mathbb{K}_W$.
When $F$ only contains \mbox{$3$--cycles}, {an $F$--factorization} with subsystems is nothing but a (nearly) Kirkman triple system with subsystems whose existence has been completely proven in \cite{PT3,PT2,PT1}. Very little is known when $F$ is any other $2$-regular graph; some recent results can be found in \cite{Prep}. More generally, we consider the definition of subsystem also for the solution to the Generalized $OP(F)$. 

In this paper we will consider both the existence of solutions to the Generalized $OP(F)$ in the case of an infinite graph $F$ that is locally finite (i.e. its vertices have finite degree) and the existence of subsystems relative to a given infinite subgraph $L$. In fact, as done by other authors for resolvable designs (see Danziger, Horsley and Webb \cite{P5}) and for $1$-factorizations (see Bonvicini and Mazzuoccolo \cite{P1}), it seems quite natural to consider the infinite case also for the Oberwolfach problem. 

In particular, in section $2$, we will recall some known facts about graphs, difference graphs and decompositions that are regular under the action of an additive group $G$. Using those results, in section $3$, we will prove that the Generalized $OP(F)$ has a solution for every infinite, locally finite graph $F$. In that proof we will use the axiom of choice: more precisely we will assume the ZFC axiomatic system. Moreover, if we consider a group $G$ that has no involutions (we will say that $G$ is involution free) of the same cardinality as $V(F)$, there exists a solution to the Generalized $OP(F)$ that is regular under the action of $G$. Finally, in section $4$, we give a characterization of the infinite subgraphs $L$ of $F$ such that there exists a solution to the $OP(F)$ that contains a subsystem relative to $L$. 
\section{Difference Graphs}
We first recall some basic definitions about graph decompositions. Let $\Lambda$ and $\Gamma$ be graphs (not necessarily finite). A $\Gamma$-decomposition of $\Lambda$ is a set $\mathcal{C}$ of copies of $\Gamma$ such that each edge of $\Lambda$ belongs to exactly one graph of $\mathcal{C}$. We say that $\mathcal{C}$ is $G$-regular if it admits $G$ as automorphism group that acts sharply transitively over the vertex set. Following \cite{ColbornDinitz07} and \cite{P04}, from now on all groups will be written additively. Given a subset $S$ of a group $G$ and an element $g\in G$ we denote by $S+g$ the set $\{x+g:\ x\in S\}$. Then we recall the following result.
\begin{prop}[See \cite{P04}]\label{righttranslations} A $\Gamma$-decomposition $\mathcal{C}$ of $\Lambda$ is $G$-regular if and only if, up to isomorphism, the following conditions hold:
\begin{itemize}
\item the vertex set of $\Lambda$ is $G$;
\item for all $B\in \mathcal{C}$ and $g\in G$, $B+g\in \mathcal{C}$. 
\end{itemize}
\end{prop} 
Clearly to describe a $G$-regular $\Gamma$-decomposition it is sufficient to exhibit a complete system $\mathcal{B}$ of representatives for the orbits of $\mathcal{C}$ under the action of $G$. The elements of $\mathcal{B}$ will be called base graphs of $\mathcal{C}$.
Here we are interested in the case where $\Lambda$ is the complete graph whose vertex set is an infinite group $G$. However the problem of finding regular decompositions has attracted much attention also in the finite case and one of the most efficient tools applied for solving this problem is the difference method (see, for example, \cite{P6b,R1,R2,R3}). 
\begin{defin}
Given a graph $\Gamma$ with vertices in an additive group $G$, the list of differences of $\Gamma$ is the multiset $\Delta(\Gamma)$ of all differences $b-b'$ between two adjacent vertices of $\Gamma$. 
\end{defin}
More generally, given a set $\mathcal{B}$ of copies of $\Gamma$ with vertices in $G$, $\Delta\mathcal{B}$ refers to the union (counting multiplicities) of all multisets $\Delta \Gamma'$, where $\Gamma'\in \mathcal{B}$.
\begin{thm}
Let $\mathcal{B}$ be a set of copies of $\Gamma$ with vertices in the group $G$.
\begin{itemize}
\item If $\Delta \mathcal{B}=G -  \{0\}$ then $\mathcal{C}=\{B+g:\ B\in\mathcal{B},g\in G\}$ is a $\Gamma$-decomposition of the complete graph $\mathbb{K}_G$ that is $G$-regular and $\mathcal{B}$ is a set of base graphs of $\mathcal{C}$.
\end{itemize}
In the case the set $\mathcal{B}$ is given by a single graph $\Gamma$ and $\Delta(\Gamma)=G -  \{0\}$, the graph $\Gamma$ is said to be a difference graph. We also denote by partial difference graph a graph $\Gamma$ such that $\Delta(\Gamma)$ is a set. 
\end{thm}
If we have a difference graph that is also a spanning subgraph of $\mathbb{K}_G$ we can say even more:
\begin{prop}\label{difference}
Let $G$ be an additive group and let $F$ be a spanning subgraph of $\mathbb{K}_G$.
If $F$ is also a difference graph then it is a base graph of a $G$-regular $F$-decomposition of $\mathbb{K}_G$ whose parts are spanning subgraphs isomorphic to $F$. In particular the Generalized $OP(F)$ has a solution.
\end{prop}
We remark that, if $F$ is a finite graph without isolated vertices, the previous proposition can not be used. In fact:
\begin{rem} Let $G$ be a finite group. Then any spanning subgraph $F$ of $\mathbb{K}_G$ that is without isolated points is not a difference graph.
\end{rem}
\proof
Indeed, if $F$ were a difference graph, we would have $\Delta(F)=G -  \{0\}$ and hence:
$$2|E(F)|=|G|-1.$$
Instead, if $F$ were a spanning subgraph of $\mathbb{K}_G$ whose vertices have finite nonzero degree:
$$2|E(F)|\geq |G|>|G|-1.$$
\endproof
In the finite case, however, it is worth recalling the similar concept of a $2$-starter (see \cite{BDT,P4,P4b,P4c,PT2b,P7}) that has strongly inspired our work.

In the infinite case, instead, the idea of Proposition \ref{difference} can be used for every nontrivial locally finite graph and indeed the main part of this paper will be dedicated to prove the following proposition:
\begin{prop}\label{main}
Let $G$ be an infinite, involution free group and let $F$ be a locally finite graph such that $|V(F)|=|E(F)|=|G|$. Then there exists an embedding $\sigma$ of $F$ in $\mathbb{K}_G$ whose image $\Gamma$ is both:
\begin{itemize}
\item a spanning subgraph $\mathbb{K}_{G}$;
\item a difference graph.
\end{itemize}
\end{prop}
Here an embedding of a graph $F$ in $\Lambda$ is defined as an injection $\sigma: V(F)\rightarrow V(\Lambda)$ such that $\{\sigma(v),\sigma(w)\}\in E(\Lambda)$ whenever $\{v,w\}\in E(F)$. The image $\sigma(F)=\Gamma$ is defined as the subgraph of $\Lambda$ such that $V(\Gamma)=\sigma(V(F))$ and $E(\Gamma)=\sigma(E(F))$.
\section{Proof of Proposition \ref{main}}
In this section we will give a proof of Proposition \ref{main}. Since we will need to apply the so called well-ordering theorem, we first recall some facts of logic.
\begin{defin}
A well-order $\prec$ on a set $X$ is a total order on $X$ with the property that every non-empty subset of $X$ has a least element.
\end{defin}
The well-order property is used for stating the following theorem that is equivalent to the axiom of choice.
\begin{thm}[Well-ordering theorem]
Every set $X$ admits a well-order $\prec$.
\end{thm}
Given an element $x\in X$ we define the section associated to it:
$$X_{\prec x}=\{y\in X: y\prec x\}.$$
\begin{cor}\label{goodgoodorder}
Every set $X$ admits a well-order $\prec$ such that the cardinality of any section is smaller than the one of $X$.
\end{cor}
\proof
Let us consider a well-order $\prec$ on $X$. Let $x$ be the smallest element such that $X_{\prec x}$ has the same cardinality as $X$. The set $Y=X_{\prec x}$ is such that all its sections with respect to the order $\prec$ have smaller cardinality. Since $Y$ instead has the same cardinality as $X$, the order $\prec$ on $Y$ induces an order $\prec'$ on $X$ with the required property.
\endproof
The well-orderings are useful because they allow proofs by induction:
\begin{thm}[Transfinite induction]\label{transfinite}
Let $X$ be a set with a well-order $\prec$ and let $P_x$ be a proposition.
Suppose that, for all $x\in X$, we have:
$$P_y| y\in X_{\prec x} \implies P_x.$$
Then $P_x$ is true for every $x\in X$.
\end{thm}
\subsection{The abelian case}
In order to better explain the proof of Proposition \ref{main} we first consider the case of abelian involution free groups. The proof for the general case will follow a similar outline but it is slightly more complicated. We begin with a technical lemma.
\begin{lem}\label{existencez}
Let $S$ be a finite subset of an involution free, abelian group $G$. Then the set of $x\in G$ such that $|\pm\{x-y: y \in S\}|<2|S|$ is finite.
\end{lem}
\proof
Given distinct $y_1$ and $y_2\in G$, for every $x\in G$, we have that $x-y_1\not=x-y_2$. Therefore the set $\pm\{x-y: y\in S\}$ has cardinality smaller than $2|S|$ if and only if $x-y_1=-(x-y_2)$ for some $y_1,y_2\in S$. Let us suppose there exist two values $x_1,x_2\in G$ such that $x_1-y_1=-(x_1-y_2)$ and $x_2-y_1=-(x_2-y_2)$. It would follow that
$$2x_1=y_1+y_2=2x_2$$ and, because of the abelianity of $G$, that $2(x_1-x_2)=0$. The latter equation is absurd since there are no involutions in $G$.
Therefore, for each pair of distinct elements $y_1,y_2\in S$, we have at most one $x$ such that $x-y_1=-(x-y_2)$. The claim follows since, set $k=|S|$, the number of such pairs is $k \choose 2$.
\endproof
Now we state a lemma that allows us to construct, inductively, an embedding $\sigma$ of the given locally finite graph $F$. First we need to provide some definitions.
In this context, we say that two vertices $p_1$ and $p_2$ of $V(F)$ are connected if there exists a path between them, and we call connected component of $v$ the set of all vertices of $V(F)$ that are connected with $v$. We define the distance between two connected vertices $p_1$ and $p_2$ as the length of the shortest path between them. Given a vertex $v\in V(F)$ the set of vertices that are at distance $1$ from $v$ will be denoted by $N(v)$ or the neighborhood of $v$; similarly, considering a subgraph $F'$ of $F$ and a vertex $v\in V(F)$ we set $N(v,F')=N(v)\cap V(F')$.
An induced subgraph $F'$ of $F$ is a subgraph obtained by vertex deletions only: $F'$ can be seen as the restriction of $F$ on the subset $W=V(F')$ of $V(F)$. We also speak about the subgraph induced by $W$. 
\begin{lem}\label{extension}
Let $G$ be an infinite, involution free, abelian group and let $F$ be a locally finite graph such that $|V(F)|=|E(F)|=|G|$. We consider an induced subgraph $F'$ of $F$ of size $|V(F')|<|V(F)|$ and let us suppose there exists an embedding $\sigma'$ of $F'$ in $\mathbb{K}_G$ as a partial difference graph $\Gamma'$. Then the following statements hold.
\begin{itemize}
\item[1)] Given $v\in V(F)$, there exists an induced subgraph $F''$ of $F$, that is an extension of $F'$ with $|V(F'')|<|V(F)|$, and an embedding $\sigma''$ of $F''$ in $\mathbb{K}_G$ as a partial difference graph $\Gamma''$ such that $\sigma''|_{F'}=\sigma'$ and $v\in V(F'')$. 
\item[2)] Given $g\in G$, there exists an induced subgraph $F''$ of $F$, that is an extension of $F'$ with $|V(F'')|<|V(F)|$, and an embedding $\sigma''$ of $F''$ in $\mathbb{K}_G$ as a partial difference graph $\Gamma''$ such that $\sigma''|_{F'}=\sigma'$ and $g\in V(\Gamma'')$. 
\item[3)] Given $g\in G -  \{0\}$, there exists an induced subgraph $F''$ of $F$, that is an extension of $F'$ with $|V(F'')|<|V(F)|$, and an embedding $\sigma''$ of $F''$ in $\mathbb{K}_G$ as a partial difference graph $\Gamma''$ such that $\sigma''|_{F'}=\sigma'$ and $g\in \Delta(\Gamma'')$. 
\end{itemize}
\end{lem}
\proof
\begin{itemize}
\item[1)] We can assume $v\in V(F) -   V(F')$. Denoted by $F''$ the graph induced by $V(F')\cup \{v\}$, we would like to extend $\sigma'$ to an embedding of $F''$ in $\mathbb{K}_G$. We note that:
\begin{itemize}
\item[$\bullet$] $|V(\Gamma')|=|V(F')|<|V(F)|=|G|$;
\item[$\bullet$] given $\delta\in \Delta(\Gamma)$ and $w\in N(v,F')$, there exist at most two values of $x$ such that one of the differences $x-\sigma'(w)$ and $\sigma'(w)-x$ is $\delta$. Since $| \Delta(\Gamma')|< |G|$ it means that the set of $x$ such that $x-\sigma'(w)$ or $\sigma'(w)-x$ is in $\Delta(\Gamma')$, for some $w\in N(v,F')$, has size at most $2|\Delta(\Gamma')||N(v,F')|<|G|$;
\item[$\bullet$] because of Lemma \ref{existencez}, the set of $x\in G$ such that $\pm\{x-\sigma'(w): w\in N(v,F')\}$ has size smaller than $2|N(v,F')|$ is finite.
\end{itemize}
Therefore it is possible to choose a vertex $x\in G$ such that:
\begin{itemize}
\item[(a)] $x\not\in V(\Gamma')$;
\item[(b)] for every $w\in N(v,F')$ the differences $x-\sigma'(w)$ and $\sigma'(w)-x$ are not in $\Delta(\Gamma')$;
\item[(c)] the list $\pm\{x-\sigma'(w): w \in N(v,F')\}$
has size $2|N(v,F')|$.
\end{itemize}
Now we define $\Gamma''$ so that $V(\Gamma'')=V(\Gamma')\cup \{x\}$ and  $E(\Gamma'')=E(\Gamma')\cup \{\{x,\sigma'(w)\}: w\in N(v,F')\}$.
Because of the properties $(b)$ and $(c)$, $\Gamma''$ is a partial difference graph. Due to property $(a)$, $x\not\in V(\Gamma')$ and hence we can define $\sigma''(v)=x$ and $\sigma''|_{F'}=\sigma'$. We obtain an extension of $\sigma'$ that embeds $F''$ in $\mathbb{K}_G$ as the partial difference graph $\Gamma''$.
\item[2)] We can assume $g\in G -   V(\Gamma')$. We note that $|V(F')|<|V(F)|$, and, since the degree of each vertex of $F'$ in $F$ is finite, also the set of vertices $v\in F$ that have distance at most one with $V(F')$ has size smaller than that of $F$.
Therefore there exists a vertex $v\in F$ such that $N(v,F)$ is disjoint from $V(F')$. Let us denote by $F''$ the graph induced by $V(F')\cup \{v\}$. We would like to extend $\sigma'$ to an embedding of $F''$ in $\mathbb{K}_G$. Since the degree of $v$ in $F''$ is zero and $g\not\in V(\Gamma')$ we can define $\sigma''(v)=g$ and $\sigma''|_{F'}=\sigma'$. Thus we obtain an extension of $\sigma'$ that embeds $F''$ in $\mathbb{K}_G$ as the partial difference graph $\Gamma''$ given by $\Gamma'$ joined with the isolated vertex $g$.
\item[3)] We can assume $g\in G -   \Delta(\Gamma')$.
We note that $|V(F')|<|V(F)|$, and, since the degree of each vertex of $F'$ in $F$ is finite, also the set $W$ of vertices $v\in F$ that have distance at most two with $V(F')$ has size smaller than that of $F$. Calling $\bar{F'}$ the graph induced by $W$, we note that $|E(\bar{F}')|<|E(F)|=|G|$. Therefore there exists $\{v,w\}\in E(F) -   E(\bar{F}')$ with $v\in V(F) -   W$ and $w\in N(v,F)$. Since $v$ has distance at least three with $V(F')$, we have that $N(v,F)$ and $N(w,F)$ are disjoint from $V(F')$. Let us denote by $F''$ the graph induced by $V(F')\cup \{v,w\}$. 
We note that, for every $x\in V(\Gamma')$, there exist at most two $y$ such that $x-y=g$ or $-(x-y)=g$. Therefore, since $|V(\Gamma')|=|V(F')|<|V(F)|=|G|$, there exist $x$ and $y$ such that $x-y=g$ and both $x,y$ are not in  $V(\Gamma')$. We would like to extend $\sigma'$ to an embedding of $F''$ in $\mathbb{K}_G$.  Since the degree of $v$ and $w$ in $F''$ is one we can define $\sigma''(v)=x$, $\sigma''(w)=y$ and $\sigma''|_{F'}=\sigma'$. Thus we obtain an extension of $\sigma'$ that embeds $F''$ in $\mathbb{K}_G$ as the partial difference graph $\Gamma''$  given by $\Gamma'$ joined with the isolated edge $\{x,y\}$.
\end{itemize}
\endproof
Now we are able to prove Proposition \ref{main} for abelian groups.
\begin{prop}\label{mainabelian}
Let $G$ be an infinite, involution free, abelian group and let $F$ be a locally finite graph such that $|V(F)|=|E(F)|=|G|$. Then there exists an embedding $\sigma$ of $F$ in $\mathbb{K}_G$ whose image $\Gamma$ is both:
\begin{itemize}
\item a spanning subgraph of $\mathbb{K}_{G}$;
\item a difference graph.
\end{itemize}
\end{prop}
\proof
We consider a well-order $\prec$ on $G$ that satisfies the condition of Corollary \ref{goodgoodorder}.
We note that, since $|V(F)|=|G|$ we can assume $V(F)=\{v_g: g\in G\}$.

Now we prove, using the transfinite induction (see Theorem \ref{transfinite}), that there exists a family of induced subgraphs $\{F_g: g\in G\}$ of $F$ and a family of embeddings $\{\sigma_g: g\in G\}$ of $F_g\rightarrow \Gamma_g$ in $\mathbb{K}_G$, such that:
\begin{itemize}
\item $\Delta(\Gamma_g)$ is a set (i.e. $\Gamma_g$ is a partial difference graph);
\item $v_g\in V(F_g)$, $g\in V(\Gamma_g)$ and, if $g$ is nonzero, $g\in \Delta(\Gamma_g)$;
\item whenever $h\prec g$ we have that $F_h\subseteq F_g$, $\Gamma_h\subseteq \Gamma_g$ and $\sigma_g|_{F_h}=\sigma_h$.
\end{itemize}
Let us assume there exists $F_h$ and $\sigma_h$ with such properties for all $h\prec g$.
We define $F_{\prec g}$ and $\Gamma_{\prec g}$ to be $\cup_{h\prec g} F_h$ and $\cup_{h\prec g} \Gamma_h$. It is easy to see that we can define the map $\sigma_{\prec g}: F_{\prec g}\rightarrow \Gamma_{\prec g}$ so that $\sigma_{\prec g}|_{F_h}=\sigma_h$ whenever $h\prec g$. Because of Corollary \ref{goodgoodorder} the section $G_{\prec g}=\{h\in G: h\prec g\}$ has a smaller cardinality than that of $G$ and thus also the size of $F_{\prec g}$ is smaller than that of $G$. Therefore we can extend the graph $F_{\prec g}$ by applying in sequence the points $1, 2$ and $3$ of Lemma \ref{extension}. We obtain an extension $F_g$ of $F_{\prec g}$ that satisfies the required conditions.

Let us now consider the graph $\Gamma=\bigcup_{g\in G} \Gamma_g$. Because of the construction, $\Gamma$ is isomorphic to $F$. Moreover, since each $\Gamma_g$ is such that $\Delta(\Gamma_g)$ is a set, $\Delta(\Gamma)$ is also a set. Finally, we note that, for every $g\in G$, $g\in V(\Gamma_g)\subseteq V(\Gamma)$ and, if $g$ is nonzero, $g\in \Delta(\Gamma_g)\subseteq\Delta(\Gamma)$. Therefore $\Gamma$ is a copy of $F$ that is a spanning subgraph of $\mathbb{K}_{G}$ such that $\Delta(\Gamma)=G -  \{0\}$.    
\endproof
\subsection{The general case}
Now we provide the proof for a generic involution free group $G$. In that case, however, Lemma \ref{existencez} is false; in fact, as noted by an anonymous referee, there exist groups which satisfy the following property:
\begin{rem}\label{contro}
There exists an infinite, involution free group $G$ such that the map $f(x)=x+x=2x$ has fibers of infinite size. 
\end{rem}
\proof
Let $H$ be the free product of countable many involutions. Namely a set of generators for $H$ is $I=\{h_i: i\in \mathbb{N}\}$ where the elements of $I$ satisfy the relations $h_i+h_i=id$ for every $i\in \mathbb{N}$; although $H$ is not abelian, we use the additive notation in accordance with \cite{ColbornDinitz07} and \cite{P4}. Now we consider the map $\phi$ from $H \times \mathbb{Z}$ to $\mathbb{Z}_2$ defined by:
$$\phi(h,k)=\phi(h_{i_1}+h_{i_2}+\cdots+ h_{i_t},k)=(t+k)\pmod{2}$$
where $h_{i_1}+h_{i_2}+\cdots +h_{i_t}$ is the unique minimal way of writing $h$ as a sum of elements of $I$. It is easy to note that $\phi$ is a group homomorphis and let us consider the kernel $G$ of this homomorphism.
To see that $G$ is involution free, assume that there is an involution $(h,k)\in G$. Clearly $k=0$, $h=h_{i_1}+h_{i_2}+\cdots+h_{i_t}$ is an involution and $t$ is even. Since $h$ is an involution, we have that $(h_{i_1}+h_{i_2}+\cdots+h_{i_t})(h_{i_1}+h_{i_2}+\cdots+h_{i_t})=id$ which means $h_{i_1}=h_{i_{t}}, h_{i_2}=h_{i_{t-1}},\dots, h_{i_{t/2}}=h_{i_{t/2+1}}$. But then $h$ would be the identity of $H$. It follows that $G$ is involution free.
Lastly, we have that $f(h_i,1)=(h_i,1)+(h_i,1)=(id,2)$ for every $h_i\in I$; therefore there exist infinitely many elements of $G$ that are mapped to the same element by $f$.
\endproof
Hence, given the group $G$ of Remark \ref{contro}, if we consider $y_1$ and $y_2$ distinct elements of $G$ such that $y_1+y_2=(id,2)$, the equality $x-y_1=-(x-y_2)$ holds for infinitely many $x\in G$, which is in contradiction with Lemma \ref{existencez}. However, we are able to prove the following weaker version of this lemma also for nonabelian groups. Given a finite subset $S$ of $G$, we define the set 
$$V_{S}=\{x\in G: |\pm\{x-y| y\in S\}|=2|S|\}.$$
\begin{lem}\label{existencez1}
Let us consider an infinite, involution free group $G$.
\begin{itemize}
\item[1)] Given $y_1\not= y_2\in G$, then we have that $|V_{\{y_1,y_2\}}|=|G|$.
\item[2)] Let $S$ be a finite subset of $G$ such that $V_{S}$ has the same cardinality as $G$. Then the set of $z\in G$ such that $V_{S\cup \{z\}}$ has cardinality smaller than that of $G$ is finite.
\end{itemize}
\end{lem}
\proof
\begin{itemize}
\item[1)]
We note that $|\pm\{x-y_i: i\in [1,2]\}|=4$ if and only if:
\begin{itemize}
\item[$\bullet$] $y_1\not=y_2$;
\item[$\bullet$] $x-y_1\not=-(x-y_2)$ that is $2(x-y_1)\not=y_2-y_1$. 
\end{itemize}
Given $x\not\in {V}_{\{y_1,y_2\}}$, we have that $y_1-x+y_1 \in {V}_{\{y_1,y_2\}}$, in fact, since $G$ is involution free: 
$$2(y_1-x+y_1-y_1)=2(y_1-x)=-2(x-y_1)=-(y_2-y_1)\not=(y_2-y_1).$$
Since the map $x\rightarrow y_1-x+y_1$ is injective, we obtain that $|G -   V_{\{y_1,y_2\}}|\leq |V_{\{y_1,y_2\}}|$. Finally, since $G$ is infinite, it follows that $|V_{\{y_1,y_2\}}|=|G|$.
\item[2)]
We note that $x\in V_{S\cup \{z\}}$ if and only if the following conditions hold:
\begin{itemize}
\item[$\bullet$] $x\in V_{S}$;
\item[$\bullet$] $x-z\not= y-x$ for every $y\in S$.
\end{itemize}
Let us denote by $y_1,\dots,y_k$ the elements of $S$ and let us suppose, by contradiction, that there exist infinitely many distincts $z^1,\dots, z^n,\dots$ such that, for every $i\in \mathbb{N}$, we have $|V_{S\cup \{z^i\}}|<|V_{S}|=|G|$. 
Therefore, for every $i\in \mathbb{N}$, we can define a map $\pi_i$ from $V_{S} -   V_{S\cup \{z^i\}}$ to $[1,k]$ where $\pi_i(x)=j$ implies $x-z^i= y_j-x$. 

Since $|V_{S\cup \{z^i\}}|<|V_{S}|=|G|$, we have:
$$V_{S} -    \left (\bigcup_{i=1}^{k+1} V_{S\cup\{z^i\}}\right)\not=\emptyset.$$ 
Let us consider $\bar{x}\in V_{S} -    (\bigcup_{i=1}^{k+1} V_{S\cup \{z^i\}})$ and let us denote by $L$ the list \\$[\pi_1(\bar{x}),\pi_2(\bar{x}),\dots,\pi_k(\bar{x}),\pi_{k+1}(\bar{x}),\dots)]$. Since this list admits values in $[1,k]$, there exist $i',i''$ such that $\pi_{i'}(\bar{x})=\pi_{i''}(\bar{x})=j$ for some $j\in[1,k]$. This means that $\bar{x}-z^{i'}= y_j-\bar{x}=\bar{x}-z^{i''}$ but this is absurd since $z^{i'}\not=z^{i''}$.
\end{itemize}
\endproof
\begin{defin}
Let $F'$ be a subgraph of $F$ and let $\sigma': F'\rightarrow \Gamma'$ be an embedding of $F'$ into $\mathbb{K}_G$. We say that $\sigma'$ satisfies the property $\star$ if, for every $v\in V(F)$, we have $|V_{\sigma'(N(v,F'))}|=|G|$.
\end{defin}
We note that, in the case $G$ is infinite, involution free and abelian, any embedding of a subgraph $F'$ of $F$ in $\mathbb{K}_G$ as a partial difference graph satisfies the property $\star$. In fact, as a consequence of Lemma \ref{existencez}, for every $v\in V(F)$, $V_{\sigma'(N(v,F'))}$ equals $G$ without at most a finite set.
\begin{lem}\label{extensionB}
Let $G$ be an infinite, involution free group and let $F$ be a locally finite graph such that $|V(F)|=|E(F)|=|G|$. We consider an induced subgraph $F'$ of $F$ of size $|V(F')|<|V(F)|$ and let us suppose there exists an embedding $\sigma'$ of $F'$ in $\mathbb{K}_G$ as a partial difference graph $\Gamma'$ that satisfies the property $\star$, then the following hold.
\begin{itemize}
\item[1)] Given $v\in V(F)$, there exists an induced subgraph $F''$ of $F$ that is an extension of $F'$ with $|V(F'')|<|V(F)|$, and an embedding $\sigma''$ of $F''$ in $\mathbb{K}_G$ as a partial difference graph $\Gamma''$ such that $\sigma''|_{F'}=\sigma'$, $v\in V(F'')$ and $\sigma''$ satisfies the property $\star$. 
\item[2)] Given $g\in G$, there exists an induced subgraph $F''$ of $F$ that is an extension of $F'$ with $|V(F'')|<|V(F)|$, and an embedding $\sigma''$ of $F''$ in $\mathbb{K}_G$ as a partial difference graph $\Gamma''$ such that $\sigma''|_{F'}=\sigma'$, $g\in V(\Gamma'')$ and $\sigma''$ satisfies the property $\star$. 
\item[3)] Given $g\in G -   \{0\}$, there exists an induced subgraph $F''$ of $F$ that is an extension of $F'$ with $|V(F'')|<|V(F)|$, and an embedding $\sigma''$ of $F''$ in $\mathbb{K}_G$ as a partial difference graph $\Gamma''$ such that $\sigma''|_{F'}=\sigma'$, $g\in \Delta(\Gamma'')$ and $\sigma''$ satisfies the property $\star$. 
\end{itemize}
\end{lem}
\proof
\begin{itemize}
\item[1)] We can assume $v\in V(F) -   V(F')$. Denote by $F''$ the graph induced by $V(F')\cup \{v\}$, we would like to extend $\sigma'$ to an embedding of $F''$ in $\mathbb{K}_G$. We note that:
\begin{itemize}
\item[$\bullet$] $|V(\Gamma')|=|V(F')|<|V(F)|=|G|$;
\item[$\bullet$] given $\delta\in \Delta(\Gamma)$ and $w\in N(v,F')$, there exist at most two values of $x$ such that one of the differences $x-\sigma'(w)$ and $\sigma'(w)-x$ is $\delta$. Since $| \Delta(\Gamma')|< |G|$ it means that the set of $x$ such that $x-\sigma'(w)$ or $\sigma'(w)-x$ is in $\Delta(\Gamma')$, for some $w\in N(v,F')$, has size at most $2|\Delta(\Gamma')||N(v,F')|<|G|$;
\item[$\bullet$] let us consider $w\in N(v,F) -   V(F')$. Because of Lemma \ref{existencez1}(2), the set of $x$ such that $V_{\sigma'(N(w,F'))\cup\{x\}}$ has cardinality smaller than that of $G$ is finite;
\item[$\bullet$] since $\sigma'$ satisfies $\star$, the set $V_{\sigma'(N(v,F'))}$ has the same cardinality as $G$.
\end{itemize}
Therefore it is possible to choose a vertex $x\in G$ such that:
\begin{itemize}
\item[(a)] $x\not\in V(\Gamma')$;
\item[(b)] for every $w\in N(v,F')$, the differences $x-\sigma'(w)$ and $\sigma'(w)-x$ are not in $\Delta(\Gamma')$;
\item[(c)] for every $w\in N(v,F) -   V(F')$, the set $V_{\sigma'(N(w,F'))\cup\{x\}}$ has the same cardinality as $G$;
\item[(d)] $x\in V_{\sigma'(N(v,F'))}$ and hence the list $\pm\{x-\sigma'(w): w \in N(v,F')\}$
has size $2|N(v,F')|$.
\end{itemize}
Now we define $\Gamma''$ so that $V(\Gamma'')=V(\Gamma')\cup \{x\}$ and  $E(\Gamma'')=E(\Gamma')\cup \{\{x,\sigma'(w)\}: w\in N(v,F')\}$.
Because of the properties $(b)$ and $(d)$, $\Gamma''$ is a partial difference graph.
Due to property $(a)$, $x\not\in V(\Gamma')$ and hence we can define $\sigma''(v)=x$ and $\sigma''|_{F'}=\sigma'$. We obtain an extension of $\sigma'$ that embeds $F''$ in $\mathbb{K}_G$ as the partial difference graph $\Gamma''$. Moreover, because of the property $(c)$, the map $\sigma''$ satisfies the property $\star$.
\item[2)] Here we proceed as in Lemma \ref{extension}(2). In order to guarantee that the extension still satisfies $\star$, it is enough to choose the vertex $v$ of the proof of Lemma \ref{extension}(2) so that $v\in V(F)$ is at distance at least three with $V(F')$.
\item[3)] Also here we proceed as in Lemma \ref{extension}(3). In order to guarantee that the extension still satisfies $\star$, it is enough to choose the edge $\{v,w\}$ of the proof of Lemma \ref{extension}(3) so that $\{v,w\}\in E(F)$ and both $v$ and $w$ have distance at least three with $V(F')$. Then property $\star$ easily follows because, for Lemma \ref{existencez1}(1), $|V_{\{\sigma''(v),\sigma''(w)\}}|=|G|$.
\end{itemize}
\endproof
We remark that, since the property $\star$ is always satisfied by abelian groups, Lemma \ref{extensionB} is actually a generalization of Lemma \ref{extension}.

We can now prove Proposition \ref{main}, for every infinite, involution free group, in the same way as Proposition \ref{mainabelian} by using Lemma \ref{extensionB} instead of Lemma \ref{extension}.
\subsection{Solution to the infinite Oberwolfach problem}
By applying Propositions \ref{difference} and \ref{main} we are able to solve the Generalized Oberwolfach problem for every infinite, locally finite graph:
\begin{thm}\label{conc}Let $G$ be an infinite, involution free group and let $F$ be a nontrivial, locally finite, graph of the same cardinality (i.e. $|V(F)|=|G|$). Then there exists a $G$-regular solution to the Generalized $OP(F)$.
\end{thm}
\proof
Since $F$ is a locally finite, infinite graph, we have that $|V(F)|\geq |E(F)|$. We divide the proof in two cases:

CASE 1: $|E(F)|=|V(F)|=|G|$. In this case the proof immediately follows from Propositions \ref{main} and \ref{difference}.

CASE 2: $|E(F)|<|V(F)|=|G|$. Then there are at most $2|E(F)|$ vertices with nonzero degree. This means that there is a family $\mathcal{B}$ of graphs such that:
\begin{itemize}
\item any $F'\in \mathcal{B}$ is a copy of $F$ whose vertex set is $V(F)$;
\item for every $v\in V(F)$ there is at most one $F'\in \mathcal{B}$ such that the degree of $v\in F'$ is nonzero;
\item $|\mathcal{B}|=|V(F)|=|G|$.
\end{itemize}
Denoting by $\bar{F}=\bigcup_{F'\in \mathcal{B}} F'$ we have that $|E(\bar{F})|=|V(\bar{F})|=|V(F)|=|G|$ and therefore, because of CASE 1, there is a solution to the Generalized $OP(\bar{F})$ that is $G$-regular. Since $\bar{F}$ is the union of copies of $F$ and each $F'\in \mathcal{B}$ is a spanning subgraph of $\bar{F}$, we obtain a $G$-regular solution also to the Generalized $OP(F)$.
\endproof
We remark that the difference method could work also for some graphs that admit vertices with infinite degree. In particular, as noted by an anonymous referee, our proof can be easily adapted to graphs $F$ whose vertices degrees are smaller than the cardinality of $G$.

However we can not hope that this procedure works for every graph $F$ with $|V(F)|=|G|$. For example, if $F$ is a complete graph, we can embed it only as $\Gamma=\mathbb{K}_G$ and $\Delta(\mathbb{K}_G)$ is not a set. Therefore we leave open the following question. 
\begin{prob}
Given an infinite, involution free group $G$, characterize the graphs $F$ that admit a $G$-regular solution to the Generalized $OP(F)$.
\end{prob}
As a consequence of Theorem \ref{conc}, we obtain a complete solution to the infinite Oberwolfach problem:
\begin{thm}\label{conc1}Let $G$ be an infinite, involution free group and let $F$ be a $2$-regular graph of the same cardinality. Then there exists a $G$-regular solution to the $OP(F)$.
\end{thm}
In the previous theorems we require that the group $G$ has no involutions, therefore one can wonder what happens if $G$ admits a nontrivial involution. In that case Theorems \ref{conc} and \ref{conc1} would be false; in fact we have the following non existence result.
\begin{rem}
Let $G$ be an infinite group with a nontrivial involution $i$ and let $F$ be a $2$-regular graph that is the disjoint union of infinitely many odd cycles. Then there does not exist any $G$-regular solution to the $OP(F)$.
\end{rem}
\proof
Let us suppose, by contradiction, that there exists a $G$-regular solution to the $OP(F)$. Let $i$ be a nontrivial involution, let $v\in V(F)$ and let $w=v^i$ where by $v^i$ we denote the action of $i$ on the vertex $v$. We note that the edge $e=\{v,w\}$ is fixed by the action of $i$, therefore the cycle $C=(w,c_2,c_3,\dots,c_{k-1},v)$ of the decomposition that contains $e$ is fixed too. This means that $i$ does not move the vertex $c_{(k-1)/2}$. But this is absurd because $c_{(k-1)/2}$ is fixed also by $0$ and therefore the action of $G$ is not sharply vertex transitive on $V(F)$.
\endproof
More generally, given a group $G$ that has nontrivial involutions, it makes sense to consider the following problem:
\begin{prob}
Given an infinite group $G$ with nontrivial involutions, characterize the graphs $F$ that admit a $G$-regular solution to the Generalized $OP(F)$.
\end{prob}
\section{Existence of Subsystems}
Here we consider infinite graphs $F$ and $L$ that are locally finite and without isolated vertices. 
Our goal will be to determine the conditions under which there exists a solution $\mathfrak{F}$ to the $OP(F)$ over the complete graph $\mathbb{K}_V$ that contains a solution $\mathfrak{L}$ to the $OP(L)$ over $\mathbb{K}_W$ where $W$ is a subset of $V$.
Clearly we must have that $|V(L)|\leq |V(F)|$ and, since the decomposition $\mathfrak{L}$ is the restriction of $\mathfrak{F}$ on $\mathbb{K}_W$, we can assume $L$ is an induced subgraph of $F$. Then, in the case $|V(L)|<|V(F)|$, we obtain a solution to the $OP(F)$ with the required subsystem using Proposition \ref{main} and Lemma \ref{extensionB}. This procedure does not work when $|V(L)|=|V(F)|$ and thus we provide another construction for this case. We begin by generalizing Lemma \ref{existencez1}.
\begin{lem}\label{existencez2}
Let us consider an infinite, involution free group $G$ and let $H$ be a normal subgroup of $G$ of the same cardinality.
\begin{itemize}
\item[1)] Given $y_1\not= y_2\in G$, there are at least $t=\left\lceil i(G:H)/2\right\rceil$ cosets $H_1,\dots,H_t$ of $H$ in $G$ such that $|V_{\{y_1,y_2\}}\cap H_j|=|G|$ for every $j\in [1,t]$.
\item[2)] Let $S$ be a finite subset of $G$ such that $V_{S}\cap H'$ has the same cardinality as $G$ for a given coset $H'$ of $H$ in $G$. 
The set of $z\in G$ for which $V_{S\cup \{z\}}\cap H'$ has cardinality smaller than that of $G$ is finite.
\end{itemize}
\end{lem}
\proof
\begin{itemize}
\item[1)]
We note that $|\pm\{x-y_i| i\in [1,2]\}|=4$ if and only if:
\begin{itemize}
\item[$\bullet$] $y_1\not=y_2$;
\item[$\bullet$] $x-y_1\not=-(x-y_2)$ that is $2(x-y_1)\not=y_2-y_1$. 
\end{itemize}
Given $x\not\in {V}_{\{y_1,y_2\}}$ we have than $y_1-x+y_1 \in {V}_{\{y_1,y_2\}}$, in fact, since $G$ is involution free: 
$$2(y_1-x+y_1-y_1)=2(y_1-x)=-2(x-y_1)=-(y_2-y_1)\not=(y_2-y_1).$$
The map $x\rightarrow y_1-x+y_1$ is injective and, since $H$ is a normal subgroup, it maps the coset $H'$ into another coset, say $H''$; it means that $|H' - (V_{\{y_1,y_2\}}\cap H')|\leq |(V_{\{y_1,y_2\}}\cap H'')|$. Finally, since $H$ is infinite, it follows that at least one between $V_{\{y_1,y_2\}}\cap H'$ and $V_{\{y_1,y_2\}}\cap H''$ has the same size as $H$ and hence as $G$.
\item[2)]
We note that $x\in V_{S\cup \{z\}}\cap H'$ if and only if the following conditions hold:
\begin{itemize}
\item[$\bullet$] $x\in V_{S}\cap H'$;
\item[$\bullet$] $x-z\not= y-x$ for every $y\in S$.
\end{itemize}
Let us denote by $y_1,\dots,y_k$ the elements of $S$ and let us suppose, by contradiction, that there exist infinitely many distincts $z^1,\dots, z^n,\dots$ such that, for every $i\in \mathbb{N}$, we have $|V_{S\cup\{z^i\}}\cap H'|<|V_{S}\cap H'|=|G|$. 
Therefore, for every $i\in \mathbb{N}$, we can define a map $\pi_i$ from $(V_{S}  -   V_{S\cup\{z^i\}})\cap H'$ to $[1,k]$ where $\pi_i(x)=j$ implies $x-z^i= y_j-x$. 

Since $|V_{S\cup\{z^i\}}\cap H'|<|V_{S}\cap H'|=|G|$, we have:
$$\left (V_{S} -    \left (\bigcup_{i=1}^{k+1} V_{S\cup\{z^i\}}\right)\right )\cap H'\not=\emptyset.$$ 
Let us consider $\bar{x}\in (V_{S} -    (\bigcup_{i=1}^{k+1} V_{S\cup\{z^i\}}))\cap H'$ and let us denote by $L$ the list \\$[\pi_1(\bar{x}),\pi_2(\bar{x}),\dots,\pi_k(\bar{x}),\pi_{k+1}(\bar{x}),\dots)]$. Since this list admits values in $[1,k]$, there exist $i',i''$ such that $\pi_{i'}(\bar{x})=\pi_{i''}(\bar{x})$. This means that $\bar{x}-z^{i'}= y_j-\bar{x}=\bar{x}-z^{i''}$ but this is absurd since $z^{i'}\not=z^{i''}$. Therefore we have only a finite set of $z\in G$ for which $V_{S\cup\{z\}}\cap H'$ has cardinality smaller than that of $G$.
\end{itemize}
\endproof
Given a subset $W$ of $V(F)$, we denote by $F/W$ the graph obtained from $F$ by contracting (i.e. identifying) the vertices of $W$ and removing the multiple edges  and the loops. We will also denote by $\max\deg(F)$ the maximum degree of the vertices of $F$. Now we provide a generalization of property $\star$.
\begin{defin}
Let us consider an embedding $\sigma': F'\rightarrow \Gamma'$ of $F'$ into $\mathbb{K}_G$. We say that $\sigma'$ satisfies the property $\star_1$ with respect to the subset $W$ of $V(F)$ if, for every $v\in V(F)$, denoted by $d$ the degree of $v$ in $F/W$, there exist $H_1,\dots,H_{d+1}$ cosets of $G$ in $H$ such that $|V_{\sigma(N(v,F')}\cap H_j|=|G|$ for every $j\in [1,d+1]$.
\end{defin}
Now we consider infinite graphs $F$ and $L$ that are locally finite and without isolated vertices such that $L$ is an induced subgraph of $F$, $|V(L)|=|V(F)|$ and $V(F) -   V(L)$ is either the empty graph (i.e. $F=L$) or has the same cardinality as $V(F)$.
Given an involution free group $G$, we would like to embed $F$ in $\mathbb{K}_G$ as a difference graph so that $L$ is embedded as a difference graph in $\mathbb{K}_H$ where $H$ is a normal subgroup of $G$. Clearly we need that $|G|=|V(F)|$, $|H|=|V(L)|$ and $|G -   H|=|V(F) -   V(L)|$. We also assume that $\left\lceil i(G:H)/2\right\rceil> \max\deg(F/V(L))$: this technical hypothesis is satisfied for example if $F=L$ or if $i(G:H)=\infty$. Under those assumptions we obtain the following:
\begin{lem}\label{mainabelian2}
Let $F'$ be an induced subgraph of $F$ such that $|V(F')|<|V(F)|$ and let $\sigma'$ be an embedding of $F'$ in $\mathbb{K}_G$ as a partial difference graph $\Gamma'$ such that $\sigma'(F'\cap L)=\Lambda'= \Gamma'\cap H$, $\Delta(\Gamma') -   \Delta(\Lambda')\subseteq G -   H$ and $\sigma'$ satisfies the property $\star_1$ with respect to $V(L)$. Then the following statements hold.
\begin{itemize}
\item[1)]  Given $v\in V(F)$, there exists an induced subgraph $F''$ of $F$ that is an extension of $F'$ with $|V(F'')|<|V(F)|$, and an embedding $\sigma''$ of $F''$ in $\mathbb{K}_G$ as a partial difference graph $\Gamma''$ such that $v\in V(F'')$, $\sigma''|_{F'}=\sigma'$, $\sigma''(F''\cap L)=\Lambda''= \Gamma''\cap H$, $\Delta(\Gamma'') -   \Delta(\Lambda'')\subseteq G -   H$ and $\sigma''$ satisfies the property $\star_1$ with respect to $V(L)$.
\item[2)] Given $g\in G$, there exists an induced subgraph $F''$ of $F$ that is an extension of $F'$ with $|V(F'')|<|V(F)|$, and an embedding $\sigma''$ of $F''$ in $\mathbb{K}_G$ as a partial difference graph $\Gamma''$ such that $g\in V(\Gamma'')$, $\sigma''|_{F'}=\sigma'$, $\sigma''(F''\cap L)=\Lambda''= \Gamma''\cap H$, $\Delta(\Gamma'') -   \Delta(\Lambda'')\subseteq G -   H$ and $\sigma''$ satisfies the property $\star_1$ with respect to $V(L)$.
\item[3)]  Given $g\in G -  \{0\}$, there exists an induced subgraph $F''$ of $F$ that is an extension of $F'$ with $|V(F'')|<|V(F)|$, and an embedding $\sigma''$ of $F''$ in $\mathbb{K}_G$ as a partial difference graph $\Gamma''$ such that $g\in \Delta(\Gamma'')$, $\sigma''|_{F'}=\sigma'$, $\sigma''(F''\cap L)=\Lambda''= \Gamma''\cap H$, $\Delta(\Gamma'') -   \Delta(\Lambda'')\subseteq G -   H$ and $\sigma''$ satisfies the property $\star_1$ with respect to $V(L)$.
\end{itemize}
\end{lem}
\proof
\begin{itemize}
\item[1)] We can assume $v\in V(F) -   V(F')$. Denoted by $F''$ the graph induced by $V(F')\cup \{v\}$, we would like to extend $\sigma'$ to an embedding of $F''$ in $\mathbb{K}_G$. We note that:
\begin{itemize}
\item[$\bullet$] for every coset $H'$ of $H$ in $G$ we have that $|V(\Gamma')|<|G|=|H|=|H'|$; 
\item[$\bullet$] given $\delta\in \Delta(\Gamma)$ and $w\in N(v,F')$, there exist at most two values of $x$ such that one of the differences $x-\sigma'(w)$ and $\sigma'(w)-x$ is $\delta$. Since $| \Delta(\Gamma')|< |G|$ it means that the set of $x$ such that $x-\sigma'(w)$ or $\sigma'(w)-x$ is in $\Delta(\Gamma')$, for some $w\in N(v,F')$, has size at most $2|\Delta(\Gamma')||N(v,F')|<|G|$;
\item[$\bullet$] let us consider $w\in N(v,F) -   V(F')$. Because of Lemma \ref{existencez2}(2), for every $H'$ such $V_{\sigma'(N(w,F'))}\cap H'$ has the same cardinality as $G$, the set of $x$ such that $V_{\sigma'(N(w,F'))\cup\{x\}}\cap H'$ has cardinality smaller than that of $G$ is finite;
\item[$\bullet$] since $\sigma'$ satisfies $\star_1$, set $d$ the degree of $v$ in $F/V(L)$, there exist $H_1,\dots,H_{d+1}$ such that the set $V_{\sigma'(N(v,F'))}\cap H_j$ has the same cardinality as $G$ for every $j\in [1,d+1]$.
\end{itemize}
If $v\in V(L)$ we consider $H'$ to be the subgroup $H$, otherwise we consider $H'\not=H$ to be a coset such that, for every $w\in N(v,F')$, we have $\sigma'(w)\not\in H'$ and $|V_{\sigma'(N(v,F'))}\cap H'|=|G|$; that coset exists since $\sigma'$ satisfies $\star_1$. Then it is possible to choose a vertex $x\in G$ such that:
\begin{itemize}
\item[(a)] $x\in H' -   V(\Gamma')$ and hence, for every $w\in N(v,F')$, we have $\pm(\sigma'(w)-x)\in H$ if and only if $\{v,w\}\in E(L)$;
\item[(b)] for every $w\in N(v,F')$, the differences $x-\sigma'(w)$ and $\sigma'(w)-x$ are not in $\Delta(\Gamma')$;
\item[(c)] for every $w\in N(v,F) -   V(F')$, denoted by $d'$ the degree of $w$ in $F/V(L)$, there exist $H_1,\dots,H_{d'+1}$ cosets of $H$ in $G$ such that $V_{\sigma'(N(w,F'))\cup \{x\}}\cap H_j$ has the same cardinality as $G$ for every $j\in [1,d'+1]$;
\item[(d)] $x\in V_{\sigma'(N(v,F'))}\cap H'$ and hence the list $\pm\{x-\sigma'(w): w \in N(v,F')\}$
has size $2|N(v,F')|$.

\end{itemize}
We define $\Gamma''$ so that $V(\Gamma'')=V(\Gamma')\cup \{x\}$ and $E(\Gamma'')=E(\Gamma')\cup \{\{x,\sigma(w)\}: w\in N(v,F')\}$. Because of the properties $(b)$ and $(d)$, $\Gamma''$ is a partial difference graph. Due to property $(a)$, $x\not \in V(\Gamma')$ and hence we can define $\sigma''(v)=x$ and $\sigma''|_{F'}=\sigma'$. We obtain an extension of $\sigma'$ that embeds $F''$ in $\mathbb{K}_G$ as the partial difference graph $\Gamma''$. Moreover, again because of the property $(a)$, $\sigma''(F''\cap L)=\Lambda''= \Gamma''\cap H$ and $\Delta(\Gamma'') -   \Delta(\Lambda'')\subseteq G -   H$. Lastly, because of the property $(c)$, the map $\sigma''$ satisfies the property $\star_1$.
\item[2)] Here we proceed as in Lemma \ref{extensionB}(2): we guarantee that the extension still satisfies $\star_1$ by choosing $v\in V(F)$ at distance at least three with $V(F')$. Moreover, in order to have that $\sigma''(F''\cap L)=\Lambda''= \Gamma''\cap H$ and $\Delta(\Gamma'') -   \Delta(\Lambda'')\subseteq G -   H$, it is enough to choose $v$ so that $v\in V(F) -   V(L)$ in the case $g\in G -   H$ and $v\in V(L)$ otherwise. Since $|V(L)|>|V(F')|\geq |V(F'\cap L)|$ and, if $F\not=L$, $|V(F) -   V(L)|=|V(F)|$, it is possible to do that choice for $v$.
\item[3)] Also here we proceed as in Lemma \ref{extensionB}(3): we guarantee that the extension still satisfies $\star_1$ by choosing $\{v,w\}\in E(F)$ so that both $v$ and $w$ have distance at least three with $V(F')$. Since $\left\lceil i(G:H)/2\right\rceil> \max\deg(F/V(L))$, property $\star_1$ easily follows from Lemma \ref{existencez2}(1). Moreover, in the case $g\in G -   H$, we choose the edges $\{v,w\}$ and $\{x,y\}$ of the proof of Lemma \ref{extensionB}(3) so that $\{v,w\}\in E(F) -   E(L)$ and $x,y\in G -   H$. As a consequence we have that $\sigma''(F''\cap L)=\Lambda''= \Gamma''\cap H$ and $\Delta(\Gamma'') -   \Delta(\Lambda'')\subseteq G -   H$.
In the case $g\in H$, instead, we need to take $\{v,w\}\in E(L)$ and $x,y\in H$.  Since $F$ and $L$ are locally finite graphs without isolated points we have that $|E(L)|>|E(F')|\geq |E(F'\cap L)|$ and, if $F\not=L$, $|E(F) -   E(L)|=|E(F)|$; therefore it is possible to do the required choice.
\end{itemize}
\endproof
Because of Lemma \ref{existencez}, in the case the group $G$ is abelian, the property $\star_1$ is always satisfied under the other assumptions of Lemma \ref{mainabelian2}. We have also noted that the condition $\left\lceil i(G:H)/2\right\rceil> \max\deg(F/V(L))$ is satisfied when $F=L$ or $i(G:H)=\infty$. 
In the first case, since $|G -   H|=|V(F) -   V(L)|$, we also have that $G=H$ and the condition $\star_1$ reduces to condition $\star$. This means that Lemma \ref{mainabelian2} is a generalization of Lemma \ref{extensionB}. Assuming instead that $i(G:H)=\infty$, it turns out that $|V(F) -   V(L)|=|V(F)|$ and we obtain the following statement.
\begin{prop}\label{main2}
Let $F$ and $L$ be infinite, locally finite graphs without isolated vertices such that $L$ is a nontrivial induced subgraph of $F$ and $|V(F) -   V(L)|= |V(F)|$.
Given an involution free group $G$ and a normal subgroup $H$ of $G$ such that $|G|=|V(F)|$, $|H|=|V(L)|$ and, $i(G:H)=\infty$, there exists an embedding $\sigma$ of $F$ in $\mathbb{K}_{G}$ such that:
\begin{itemize}
\item the image $\Gamma$ of $F$ is a difference graph;
\item $\Gamma$ is a spanning subgraph of $\mathbb{K}_{G}$;
\item calling $\Lambda=\sigma(L)$, we have that $\Delta(\Lambda)=H$;
\item $V(\Lambda)=H$.
\end{itemize} 
\end{prop}
\proof
We divide the proof in two cases.

CASE 1: $|V(L)|<|V(F)|$. Because of Proposition \ref{mainabelian} there exists an embedding $\sigma':L\rightarrow \Lambda$ in $\mathbb{K}_{H}$ such that $\Delta(\Lambda)=H -   \{0\}$ and $V(\Lambda)=H$. We consider a well-order $\prec$ on $G -   H$ that satisfies the condition of Corollary \ref{goodgoodorder}.
We note that, since $|V(F) -   V(L)|=|G -   H|$, we can assume $V(F) -   V(L)=\{v_g: g\in G -   H\}$.
Then, following the proof of Proposition \ref{main}, by repeated application of Lemma \ref{extensionB}, we get an ascending family $\sigma_g$ of extensions of $\sigma'$.
We obtain an embedding $\sigma$ of $F$ in $\Gamma\subseteq \mathbb{K}_{G}$ such that $\Gamma$ is a difference graph, $\Gamma$ is a spanning subgraph of $\mathbb{K}_{G}$ and, since $L$ is an induced subgraph of $F$, $\sigma|_{F\cap V(L)}=\sigma'$.

CASE 2: $|V(L)|=|V(F)|$. In this case we follow the proof of Proposition \ref{main} extending, at each step, the embedding $\sigma_g$ via Lemma \ref{mainabelian2} instead of via Lemma \ref{extension}. We obtain an embedding $\sigma$ of $F$ in $\Gamma\subseteq \mathbb{K}_{G}$ such that $\Gamma$ is both a difference graph and a spanning subgraph of $\mathbb{K}_{G}$; set $\Lambda=\sigma(L)$, since $L$ is an induced subgraph of $F$, $\Lambda$ is the spanning subgraph of $\mathbb{K}_{H}$ given by $H\cap \Gamma$ and $\Delta(\Lambda)=H -  \{0\}$.
\endproof 
Proposition \ref{main2} can be used in order to characterize the existence of subsystems of the Oberwolfach problem in the infinite case.
\begin{thm}\label{las}
Let $F$ and $L$ be two infinite, locally finite graphs without isolated vertices. Then there exists a solution to the Generalized $OP(F)$ that admits a subsystem relative to $L$ if and only if $F$ contains a copy $\tilde{L}$ of $L$ such that: 
\begin{itemize}
\item[1)] $\tilde{L}$ is an induced subgraph of $F$;
\item[2)] either $\tilde{L}=F$ or we have that $|V(F) -   V(\tilde{L})|=|V(F)|$.
\end{itemize}  
\end{thm}
\proof
Let us suppose conditions $1$ and $2$ hold. Let $G$ be an involution free group with the same cardinality as $V(F)$ and let $H$ be a normal subgroup of $G$ such that $|H|=|V(L)|$ and $i(G:H)=\infty$. Then, according to Proposition \ref{main2}, there exists an embedding $\sigma$ of $F\rightarrow \Gamma$ in $\mathbb{K}_{G}$ such that:
\begin{itemize}
\item $\Gamma$ is a difference graph;
\item $\Gamma$ is a spanning subgraph of $\mathbb{K}_{G}$;
\item set $\Lambda=\sigma(\tilde{L})$, $\Delta(\Lambda)=H -  \{0\}$;
\item $\Lambda$ is a spanning subgraph of $\mathbb{K}_{H}$.
\end{itemize}
Because of Proposition \ref{difference}, the action of $G$ over $\Gamma$ gives us a solution to the Generalized $OP(F)$ while the action of $H$ over $\Lambda$ gives us a solution to the Generalized $OP(L)$ that is contained in the first one.

Let us now suppose there exists a solution $\mathfrak{F}$ to the Generalized $OP(F)$ over the graph $\mathbb{K}_{V(F)}$, that contains a subsystem $\mathfrak{L}$ relative to $L$. This means that the vertex set of any element of $\mathfrak{L}$ is a subset $W$ of $V(F)$. We can assume that $\mathfrak{F}=\{F_{\alpha}: \alpha\in \mathcal{F}\}$ and $\mathfrak{L}=\{L_{\alpha}: \alpha\in \mathcal{L}\}$ where $\mathcal{L}$ is a subset of $\mathcal{F}$ and $L_{\alpha}=F_{\alpha}|W$. Moreover, given $\bar{\alpha}\in \mathcal{L}$, we can identify $F$ with $F_{\bar{\alpha}}$ and we define $\tilde{L}$ to be the graph $L_{\bar{\alpha}}$. Clearly, since $L_{\bar{\alpha}}$ is the restriction of $F_{\bar{\alpha}}$ on the vertex set $W=V(L_{\bar{\alpha}})$, $\tilde{L}$ is an induced subgraph of $F$.

Finally, let us suppose by contradiction, that condition $2$ does not hold for the graph $L_{\bar{\alpha}}=\tilde{L}$, and hence $|V(F) -   V(\tilde{L})|< |V(F)|$ which implies $|V(F)|=|V(\tilde{L})|$. Let us consider a vertex $v\in \tilde{L}$, in the graphs of $\{F_{\alpha}| \alpha  \in \mathcal{F} -  \mathcal{L}\}$, $v$ is connected with all (and only with) the vertices of $V(F) -   V(\tilde{L})$. Since the degree of $v$ in any graph of $\{F_{\alpha}| \alpha  \in \mathcal{F} -  \mathcal{L}\}$ is nonzero, we obtain that $|\mathcal{F} -  \mathcal{L}|\leq |V(F) -   V(\tilde{L})|$. Similarly, given a vertex $w\in V(F) -   V(\tilde{L})$, in the graphs of $\{F_{\alpha}| \alpha  \in \mathcal{F} -  \mathcal{L}\}$, $w$ is connected with all (and not necessarily only with) the vertices of $V(\tilde{L})$. 
We note that $w$ has a finite degree in each graph of $\{F_{\alpha}| \alpha  \in \mathcal{F} -  \mathcal{L}\}$ and $|V(\tilde{L})|$ is infinite; this means that $|V(\tilde{L})|\leq |\mathcal{F} -  \mathcal{L}|$.
Thus $|V(F)|=|V(\tilde{L})|\leq |\mathcal{F} -  \mathcal{L}|\leq |V(F) -   V(\tilde{L})|$. But this is absurd because we have assumed $|V(F) -   V(\tilde{L})|< |V(F)|$. Therefore the condition $2$ is necessary in order to have the required subsystem.
\endproof
The proof of Theorem \ref{las} provide us regular solutions to the Generalized Oberwolfach problem that contains regular subsystems. More precisely:
\begin{rem}
Let $F$ and $L$ be infinite, locally finite graphs without isolated vertices such that $L$ is a nontrivial induced subgraph of $F$ and $|V(F) -   V(L)|= |V(F)|$.
Given an involution free group $G$ and a normal subgroup $H$ of $G$ such that $|G|=|V(F)|$, $|H|=|V(L)|$ and, $i(G:H)=\infty$, there exists a $G$-regular solution to the Generalized $OP(F)$ that contains an $H$-regular solution to the Generalized $OP(L)$.
\end{rem}
We note that, given a $2$-factor $F$ and a subgraph $L$ of $F$ that is also a $2$-factor, then $L$ is always an induced subgraph of $F$. Therefore we obtain the following characterization to the existence of subsystems of solutions to the Oberwolfach problem.
\begin{thm}\label{Obersub}
Let $F$ and $L$ be two infinite $2$-regular graphs. Then there exists a solution to the $OP(F)$ that admits a subsystem relative to $L$ if and only if:
\begin{itemize}
\item $F$ contains a copy $\tilde{L}$ of $L$ such that either $\tilde{L}=F$ or $|V(F) -   V(\tilde{L})|=|V(F)|$.
\end{itemize}  
\end{thm}
Finally, we remark that, here, we have assumed the degree of each vertex of $F$ and of $L$ to be finite and nonzero. We have used that hypothesis, both in the proof of Lemma \ref{mainabelian2} and in the one of Theorem \ref{las}, in order to avoid isolated vertices. As a consequence we leave the following question open:
\begin{prob}
Characterize the infinite graphs $F$ and $L$ such that there exists a solution to the Generalized $OP(F)$ that admits a subsystem relative to $L$.
\end{prob}
\section*{Aknowledgments}
The author would like to thank the anonymous referees for they many valuable comments and suggestions.
Research for this paper was partially carried out while the author was visiting the Beijing Jiaotong University. He expresses his sincere thanks to the 111 Project of China [grant number B16002] for financial support and to the Department of Mathematics of the Beijing Jiaotong University for their kind hospitality. 

\end{document}